
\magnification=1200

\input amstex
\documentstyle{amsppt}

\vcorrection{-.4in}\hcorrection{.2 in}

\input xypic
\input epsf

\define\CDdashright#1#2{&\,\mathop{\dashrightarrow}\limits^{#1}_{#2}\,&}
\define\CDdashleft#1#2{&\,\mathop{\dashleftarrow}\limits^{#1}_{#2}\,&}

\def\P{\Bbb P}

\def\G{\Bbb G}

\def\Til#1{{\widetilde{#1}}}

\def\PGL{\text{PGL}}

\def\im{\text{\rm im}\,}

\def\siltable#1.{
\vbox{\tabskip=0pt \offinterlineskip
\halign to360pt{\strut##& ##\tabskip=1em plus2em&
  \hfil##\hfil& \vrule##&
  \hfil##\hfil& \vrule##&
  \hfil##\hfil& \vrule##&
  \hfil##\hfil& \vrule\thinspace\vrule##&
  \hfil##\hfil& \vrule\thinspace\vrule##&
  \hfil##\hfil& ##\tabskip=0pt\cr
#1}}}

\CompileMatrices

\topmatter
\title Plane curves with small linear orbits II\endtitle
\author Paolo Aluffi\\ Carel Faber\endauthor
\date October 1999\enddate
\address 
Mathematics Department, 
Florida State University,
Tallahassee FL 32306, U.S.A.
\endaddress
\email aluffi\@math.fsu.edu\endemail
\address
Dept.~of Mathematics,
Oklahoma State University,
Stillwater OK 74078, U.S.A.
\endaddress
\email cffaber\@math.okstate.edu\endemail
\address
Department of Mathematics,
KTH,
100 44 Stockholm, Sweden
\endaddress
\email carel\@math.kth.se\endemail
\abstract The `linear orbit' of a plane curve of degree $d$ is its orbit 
in $\P^{d(d+3)/2}$ under the natural action of $\PGL(3)$. We classify
curves with positive dimensional stabilizer, and we compute the degree
of the closure of the linear orbits of curves supported
on unions of lines. Together with the results of \cite{3},
this encompasses the enumerative geometry of all plane curves with small 
linear
orbit. This information will serve elsewhere as an ingredient in the 
computation of
the degree of the orbit closure of an arbitrary plane curve.
\endabstract


\endtopmatter

\document
\footnote[]{Both authors gratefully acknowledge partial NSF support,
under grants DMS-9500843 and DMS-9801257.}
\leftheadtext{Paolo Aluffi, Carel Faber}

\head \S0. Introduction\endhead

The `linear orbit' of a plane curve is its orbit under the natural
action of $\PGL(3)$ on $\P^N=\P^{d(d+3)/2}$, the space parameterizing
plane curves of degree $d$. We say that a plane curve has `small'
linear orbit if its stabilizer has positive dimension, so that the
dimension of the orbit is less than the expected dimension
$8=\dim\PGL(3)$. This paper is concerned with such curves: we classify
them (\S1), then compute the degree of the closure
of the orbits of such curves as well as
other enumerative information (\S2). The enumerative
computations in this paper regard curves supported on unions of lines;
if all lines but at most one are concurrent, the configuration has
small orbit. The enumerative geometry of all other curves with small
orbit (that is, of curves with small orbit and containing some
non@-linear component) is substantially more involved, and is studied
in \cite{3}.

Orbits of {\it smooth\/} plane curves are studied in \cite{2}. The
results of this paper, together with the results in \cite{3}, form
an essential ingredient towards the computation of the degree of the
orbit closure for an arbitrary plane curve \cite{5}.
For a slightly more expanded discussion of the general context, we
refer the reader to the introduction of \cite{3}.
The enumerative computations in this paper rely on the strategy
employed in  \cite{2} and \cite{3}:
we construct a smooth complete variety dominating the orbit closure,
and perform the necessary intersection theory. The variety is
obtained by suitably blowing up the $\P^8$ of $3\times 3$ matrices,
viewed as a compactification of $\PGL(3)$: the action on a given curve
extends to a rational map on this $\P^8$, and the variety is obtained
by performing a sequence of blow@-ups removing 
the indeterminacies of this map. We have taken
the {\it degree\/} of the orbit closure of a curve as the main focus of
our work, but it should be noted that explicit constructions of good
varieties dominating an orbit closure in principle allow the
computation of other invariants, such as the multiplicity of the orbit
closure along its boundary, or various kinds of characteristic classes.

As the referee pointed out, the subject of plane curves with positive
dimensional stabilizer is very classical, going back as it does
to F.~Klein and S.~Lie \cite{7}.
Their classification of `W-curves' includes transcendental curves
as well, and essentially coincides with ours in the algebraic case.
Our viewpoint is however somewhat different, as we need to compile
information regarding the degree of the stabilizers of the curves,
especially when these are reducible. We do not know of any work
on the enumerative geometry of these curves previous to ours.

{\bf Acknowledgments.\/} P.A.~thanks Florida State University for a
`Developing Scholar Award', and
C.F.~thanks the Max-Planck-Institut f\"ur Mathematik, Bonn, for 
support during the final stages of the preparation of this paper.

\def\Aut{\text{Stab}}

\head \S1. Classification of curves with small orbits\endhead

In this section we classify the plane curves with small orbits, i.e.,
the plane curves with positive dimensional stabilizer. We work
over an algebraically closed field of characteristic $0$.
See \S2.4 for a few remarks about linear orbits in positive
characteristic.

Let $C$ be a plane curve with $\dim \Aut(C)>0$. (By $\Aut(C)$ we mean 
the subgroup of $\PGL(3)$ of transformations fixing $C$.) Denote 
by $G$ the connected component of the identity of $\Aut(C)$. It is 
clear that $G$ fixes the irreducible components of $C$ and does not 
depend on the multiplicities of these components.

Every curve that is not just a line contains 4 points that form a frame; hence
$G$ acts faithfully on such curves. It follows that 
$C$ contains
infinitely many points which are not fixed by $G$. By imposing
the condition that such a point be fixed, we find a subgroup of $G$ of
dimension one less; repeating this, we find subgroups of $G$ of every
dimension between 0 and $\dim G$, so that in particular $G$ contains a
1-dimensional subgroup. It is well-known that a connected 1-dimensional
linear algebraic group is isomorphic either to the additive group $\G_a$ 
or to the multiplicative group $\G_m\,$.

First we consider the case where $G$ contains a $\G_m\,$. It can be 
diagonalized: denoting 
the embedding of $\G_m$ into $G$ by $\gamma$,
we can find
homogeneous coordinates $(x:y:z)$ on $\P^2$ such that the effect of 
$\gamma(t)$ is given by $(x:y:z)\mapsto (x:t^a y:t^b z)$, with $0\le 2a\le b$
and $a$ and $b$ coprime. It is easy to see that all irreducible curves
that are fixed by $\gamma(\G_m)$ are given by the equations $x$, $y$, $z$
or $y^b+\lambda z^a x^{b-a}$ with $\lambda\neq0$. 

It follows that $C$ consists of irreducible components of the form above, with
arbitrary multiplicities. We list here the types of curves $C$ so obtained, 
together with 
the dimension of the orbit ${\Cal O}_C$.
\roster
\item $C$ consists of a single line; $\dim{\Cal O}_C=2$ (this is the only 
case where $\Aut(C)$ doesn't act faithfully on $C$).
\item $C$ consists of 2 (distinct) lines; $\dim{\Cal O}_C=4$.
\item $C$ consists of 3 or more concurrent lines; $\dim{\Cal O}_C=5$. (We call
this configuration a {\it star\/}.)
\item $C$ is a triangle (consisting of 3 lines in general position);
$\dim{\Cal O}_C=6$.
\item $C$ consists of 3 or more concurrent lines, together with 1 other
(non-concur\-rent) line; $\dim{\Cal O}_C=7$. (We call this configuration a
{\it fan\/}.)
\item $C$ consists of a single conic; $\dim{\Cal O}_C=5$.
\item $C$ consists of a conic and a tangent line; $\dim{\Cal O}_C=6$.
\item $C$ consists of a conic and 2 (distinct) tangent lines;
$\dim{\Cal O}_C=7$.
\item $C$ consists of a conic and a transversal line and may contain either
one of the tangent lines at the 2 points of intersection or both of them;
$\dim{\Cal O}_C=7$.
\item $C$ consists of 2 or more bitangent conics (conics in the pencil
$y^2+\lambda x z$) and may contain the line $y$ through the two points
of intersection as well as the lines $x$ and/or $z$, tangent lines to the
conics at the points of intersection; again, $\dim{\Cal O}_C=7$.
\item $C$ consists of 1 or more (irreducible) curves from the pencil
$y^b+\lambda z^a x^{b-a}$, with $b\ge 3$, and may contain the lines $x$
and/or $y$ and/or $z$; $\dim{\Cal O}_C=7$.
\endroster

Next we consider the case where $G$ contains a $\G_a\,$. The $\G_a$-orbit
of a point which is not a fixed point is an affine curve. Upon taking the
projective closure of this curve, we find a fixed point of the $\G_a$-action
as well as a fixed (tangent) line through that point. Hence in suitable
coordinates $(x:y:z)$ the image of the embedding $\gamma:\G_a\to G$ 
consists of lower
triangular matrices; it is easy to see that the diagonal entries are 
equal to 1.
$$
\gamma(t)=\pmatrix
1 & 0 & 0 \\
a(t) & 1 & 0 \\
b(t) & c(t) & 1
\endpmatrix
$$
Here $a(t)$, $b(t)$ and $c(t)$ are polynomials. Necessary and sufficient
conditions for $\gamma$ to be a homomorphism of algebraic groups are that
$a$ and $c$ are additive polynomials ($a(s+t)=a(s)+a(t)$, $c(s+t)=c(s)+c(t)$)
and that $b$ satisfies $b(s)+b(t)+a(s)c(t)=b(s+t)$.

Since the characteristic is $0$, both 
$a$ and $c$
are necessarily linear, forcing $b$ to be quadratic; the general solution
can be written in the form
$$
\gamma(t)=\pmatrix
1 & 0 & 0 \\
at & 1 & 0 \\
bt+{1\over2}act^2 & ct & 1
\endpmatrix\,.
$$
Since $\gamma$ is an embedding, the constants $a$, $b$ and $c$ cannot all
vanish. When $a$ and $c$ are both non-zero, the only fixed point is $(0:0:1)$.
When $c=0$, the locus of fixed points is the line $x=0$; when $a=0$, the locus
of fixed points is the line $bx+cy=0$. It follows that every irreducible
component of $C$ has degree 1 or 2; if it doesn't consist entirely of fixed
points, it has a unique fixed point. If $C$ consists only of lines, these
necessarily form a star. If $C$ contains a unique conic, it may also contain
the tangent line at one point; again, we find no new configurations this way.
If however $C$ contains two or more conics, every conic intersects every
other conic in exactly one point (with multiplicity 4); this point,
$(0:0:1)$, lies on all the conics. The curve $C$ may also contain the
(tangent) line $x=0$. (The conics are elements of the pencil through
$cy^2+2x(by-az)$ and $x^2$.) We find thus one new configuration:
\roster
\item"(12)" $C$ contains 2 or more 
conics from a pencil through a conic and a
double tangent line; it may also contain that tangent line. In this case,
$\dim{\Cal O}_C=7$ (the identity component of $\Aut(C)$ is $\G_a$).
\endroster

The following picture represents schematically the curves described
above.
$$\epsffile{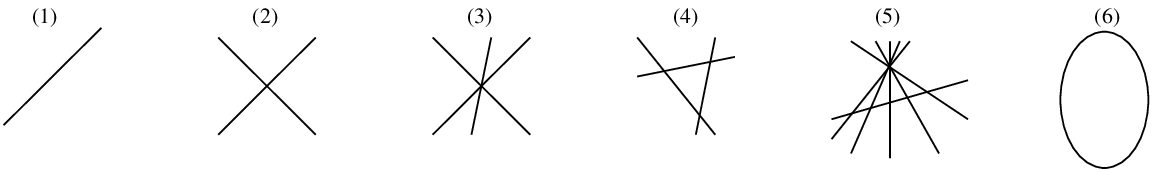}$$
$$\epsffile{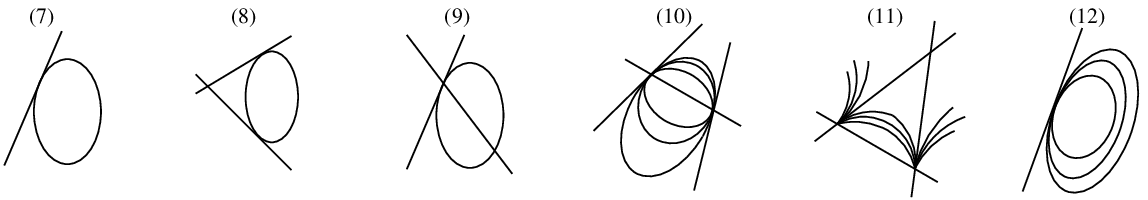}$$

We conclude this section with a description of the stabilizers
$\Aut(C)$ of the various types of plane curves $C$ with small orbit;
the degree of the stabilizer is one ingredient in the computation of
the degree of the orbit closure of the corresponding curve.

Choose coordinates $(x:y:z)$ on $\P^2$,
$(p_0 : p_1 : p_2 : p_3 : p_4 : p_5 : p_6 : p_7 : p_8)$ on $\P^8$,
and $(a:b:c:d)$ on $\P^3$.
The list below corresponds to the list above.
\roster
\item The stabilizer is 
(the intersection of $\PGL(3)$ and)
a linear subspace of dimension 6. 
For $C$ with equation $x^d$ it is given by $p_1=p_2=0$.
\item The group $G$, the connected component of the identity of
$\Aut(C)$, is a linear subspace of dimension 4. 
For $C$ with equation $x^ky^{d-k}$ it is given by
$p_1=p_2=p_3=p_5=0$. The index of $G$ in $\Aut(C)$ is $2$ when $2k=d$
and $1$ otherwise.
\item The group $G$ is a linear subspace of dimension 3. 
For $C$ with equation $f(x,y)$ (with 3 or more distinct linear factors)
it is given by $p_0=p_4$ and $p_1=p_2=p_3=p_5=0$. 
The quotient $\Aut(C)/G$ is isomorphic to the stabilizer in $\PGL(2)$
of the $d$@-tuple with equation $f(x,y)$, via the map 
sending a $3\times3$ matrix
$(p_0 : p_1 : 0 : p_3 : p_4 : 0 : p_6 : p_7 : p_8)$
in $\Aut(C)$ to the $2\times2$ matrix
$(p_0 : p_1 : p_3 : p_4)$
in the stabilizer of the $d$@-tuple.
\item The group $G$ is a linear subspace of dimension 2.
For $C$ with equation $x^ay^bz^c$ it is given by
$p_1=p_2=p_3=p_5=p_6=p_7=0$. The quotient $\Aut(C)/G$ 
is isomorphic to the symmetry group of the triple $(a,b,c)$.
\item The group $G$ is a line in $\PGL(3)$. 
For $C$ with equation $f(x,y)z^k$ it is given by
$p_0=p_4$ and $p_1=p_2=p_3=p_5=p_6=p_7=0$.
As in (3), the quotient $\Aut(C)/G$ is isomorphic to the stabilizer 
in $\PGL(2)$ of the $(d-k)$@-tuple with equation $f(x,y)$.
\item The stabilizer of a conic $x^2+yz$ is the image of the
embedding of $\PGL(2)$ into $\PGL(3)$ given by
$(a:b:c:d)\mapsto(ad+bc : bd : -ac : 2cd : d^2 : -c^2 : -2ab : -b^2 : a^2)$.
Its closure has degree 8 since it is a projection 
of the second Veronese of $\P^3$
from the $\P^9$ of
quadrics in $(a:b:c:d)$, 
from the
determinant quadric $ad-bc$.
\item The stabilizer of a conic with tangent line $y(x^2+yz)$
is the image of the embedding
$(a:b:0:d)\mapsto(ad : bd : 0 : 0 : d^2 : 0 : -2ab : -b^2 : a^2)$.
Its closure has degree 4 since it is the second Veronese of $\P^2$.
\item For a conic with two tangent lines $y^kz^l(x^2+yz)^m$, the
group $G$ consists of the diagonal matrices with entries
$(ad,d^2,a^2)$. Its closure has degree 2. The quotient $\Aut(C)/G$
is isomorphic to the symmetry group of the pair $(k,l)$.
\item\hskip-1.8mm, (10), (11) The stabilizers in these cases are 
discussed in \cite{3}, Lemma~3.1.
\item"(12)" These curves and their stabilizers are discussed
in \cite{3}, \S4.1.
\endroster

\head 2. Predegree polynomials for configurations of lines\endhead

\subheading{\S2.1 Predegree polynomials} We move now to the
enumerative geometry of these curves. Our main task is to compute
the {\it degree\/} of the orbit closure of a given curve with small
orbit. We have chosen this question because of its clear enumerative
interpretation, and because the tool needed to solve it (that is, a
nonsingular variety dominating the orbit closure) could in principle
allow us to compute several other invariants of the orbit closure. 
If the curve $C$ is reduced and its orbit has dimension $j$, then the
degree of the orbit closure of $C$ is the number of translates of $C$
which contain $j$ points in general position in the plane. In very
special cases (for example, when $C$ is a triangle) this number can be
computed by naive combinatorial considerations. In general this is not
possible. 

The computation of the degree of the orbit closure of a curve with
small orbit and
containing non@-linear components (that is, items (6) through (12) in
the classification of \S1) requires a lengthy study, paralleling the
embedded resolution of such curves. As this study is rather involved,
we have devoted to it a separate paper (\cite{3}). Here we deal with
curves consisting of configurations of (possibly multiple) lines. We
have seen in \S1 that, among these, the ones with small orbit are the
configurations in which all but at most one of the lines are
concurrent (types (1) through (5)). We call a set of concurrent
lines a {\it star\/}, and the union of a set of concurrent lines with a
general line a {\it fan\/}.

The machinery necessary to study the enumerative geometry of stars and
fans yields in fact analogous results for {\it arbitrary\/} unions of
lines (whose orbits are not necessarily small); we will state the more
general result in \S2.3, and derive from it formulas for stars and
fans.

As pointed out already in \cite{3}, a more refined degree question
is natural in our framework. In view of the applications of the
results we obtain here and in \cite{3}, it is useful to compute the
degree of the orbit closure as well as of the subsets of the orbit
closure determined by imposing linear conditions on the transformation
applied to the given curve.

To state this question precisely, let $C$ be a plane curve, with
equation $F(x:y:z)=0$. The orbit of $C$ is the image of the action map
$\PGL(3) @>>> \P^N$
sending $\varphi$ to the curve with equation $F(\varphi(x):
\varphi(y): \varphi(z))=0$. Embedding $\PGL(3)$ in the $\P^8$ of
$3\times 3$ matrices, the action extends to a rational map
$$\P^8 \overset c\to\dashrightarrow \P^N\quad,$$
and the orbit closure of $C$ is the closure of the image of this
rational map. `Imposing general linear conditions' on $\varphi\in\PGL(3)$
amounts to restricting this map $c$ to general linear subspaces of
$\P^8$.

We will compute here the degree of the closure of the image of general
subspaces of $\P^8$ of all dimensions, for $C=$ an arbitrary union of
(possibly multiple) lines. For $0\le j\le 8$, consider a general
$\P^j$ in $\P^8$, and let $f_j$, $d_j$ denote respectively the number
of points in the general fiber of $c_{|\P^j}$ and the intersection
number of $c(\P^j)$ and a codimension@-$j$ linear subspace of $\P^N$.
The products $f_j\cdot d_j$ are the {\it predegrees\/} of $C$. We
assemble these numbers into a generating function
$$\sum_{j\ge 0}\left(f_j\cdot d_j\right)\frac{t^j}{j!}\quad,$$
which we call the `adjusted predegree polynomial' of $C$.

\remark{Remarks} 
(1) If the orbit of $C$ has dimension $j$, then 

--- $d_k=0$ for $k>j$, and $d_j$ is the degree of the orbit closure;

--- $f_j$ is the degree of the stabilizer $\Aut(C)$ of $C$, and
    $f_k=1$ for $k<j$;

--- since $C$ consists of lines, the closure of 
every connected component of the stabilizer of $C$ is a 
linear subspace of $\P^8$; hence $f_j$ is the number of
connected components of the stabilizer of $C$.
For example, $f_j=1$ if the lines in the configuration appear
with distinct multiplicities.

(2) $C$ has small orbit if and only if its adjusted
predegree polynomial has degree $<8$.
By (1), the degree of the
orbit closure of $C$ can be recovered from the leading
coefficient of the adjusted predegree polynomial of $C$ 
by dividing by the
number of components of the stabilizer of $C$.

(3) The denominators are introduced in the definition of the adjusted
predegree polynomial because this uncovers a certain amount of
structure which would otherwise be lost (and is completely invisible
in the individual predegrees). This structure amounts to a
multiplicativity of the adjusted predegree polynomial with respect to
union of transversal lines, and allows a convenient shortcut in the
computation of the predegrees of stars and fans, 
cf.~Lemma 2.6 and Corollary~2.11.
\endremark

\vskip6pt
Computing adjusted predegree polynomials (rather than degrees alone)
allows us to deal at once with all configurations of lines, regardless
of the dimension of their stabilizer. Adjusted predegree polynomials
for all other small orbits are computed in \cite{3}, \S4.2.

{From} a geometric point of view, our approach to the problem will be to
construct (in \S2.2)
a nonsingular variety $\Til V$ resolving the indeterminacies of the
rational map~$c$:
$$\diagram
& {\Til V} \dlto \drto^{\tilde c}\\
{\P^8} \xdashed[0,2]^{c}|>\tip & & {\P^N} \\
\enddiagram$$

Once this is accomplished, we can let $\Til H$, $\Til W$ respectively
be the pull@-backs of the hyperplane class $H$ in $\P^8$ and $W$ in
$\P^N$;
then tracing the definitions shows that the adjusted predegree
polynomial equals
$$\sum_{j\ge 0} \left(\int \Til H^{8-j} \Til W^j\right)
\frac{t^j}{j!}\quad.$$
The construction of \S2.2 and intersection theory allow us in
\S2.3 to compute the individual intersection degrees
$\int \Til H^{8-j} \Til W^j$
of classes in $\Til V$.

\subheading{\S2.2 Blow@-ups} Let $C$ be a degree@-$d$ plane curve
consisting of a union of lines: that is, the ideal of $C$ is generated
by a product
$$F(x,y,z)=\prod L_i(x,y,z)^{r_i}\quad, \quad \sum r_i=d$$
where the $L_i$ are distinct linear homogeneous polynomials.
In this section we resolve the indeterminacies of the corresponding
rational map $c:\P^8 \dashrightarrow \P^N$, $N=d(d+3)/2$. This
construction does not use any information regarding whether the orbit
of $C$ is small or not.

First, observe that the base locus of $c$ reflects the combinatorics
of the configuration of lines. More precisely, $\varphi\in \P^8$ is in
the base locus if and only if $F(\varphi(x,y,z))$ is identically~0,
that is, if the image of $\varphi$ is contained in one of the lines of
the configuration. The condition $\im\varphi\subset$~line specifies a
certain $\P^5$ in the $\P^8$ of matrices. Hence, the base locus
consists of a union of $\P^5$'s, each corresponding to a line of the
configuration; the intersection of two $\P^5$'s corresponding to two
lines $\ell$, $m$ is the $\P^2$ of matrices whose image is the point
$\ell\cap m$.

Summarizing, the line configuration determines a set of disjoint
$\P^2$'s in $\P^8$, one for each point of intersection of two lines in
the configuration, and a set of $\P^5$'s in $\P^8$, one for each line
in the configuration. This is the base locus of $c$.

\proclaim{Proposition 2.1} Let $\Til V @>\pi>> \P^8$ be the variety
obtained by first blowing up $\P^8$ along the distinguished $\P^2$'s,
and then along the proper transforms of the distinguished $\P^5$'s.
Then $\Til V$ resolves the indeterminacies of $c$; that is, there is a
commutative diagram
$$\diagram
& {\Til V} \dlto_{\pi} \drto^{\tilde c}\\
{\P^8} \xdashed[0,2]^{c}|>\tip & & {\P^N} \\
\enddiagram$$
with $\tilde c$ a regular map.\endproclaim

The rest of this section is devoted to the proof of this proposition.

Given $C$ as above, every point $(x_0:y_0:z_0)$ in the plane
determines a hypersurface in $\P^8$, with ideal generated by
$F(\varphi(x_0,y_0,z_0))\,:$
concretely, this hypersurface consists of the matrices $\varphi$ such
that the translate of $C$ by $\varphi$ contains the point
$(x_0,y_0,z_0)$ (or is undefined). We call such hypersurfaces (and
their proper transforms) `point@-conditions'. Note that the rational
map $c$, and its lift to varieties dominating $\P^8$, is precisely the
map determined by the linear system generated by the
point@-conditions. What we have to show is that the (proper transforms
of the) point@-conditions in $\Til V$ generate a base@-point@-free
linear system.

Let $V'$ be the variety obtained after the first set of blow@-ups.
\proclaim{Claim 2.2} (i) The proper transforms of the distinguished $\P^5$'s
in $V'$ are disjoint.

(ii) The multiplicity of a point@-condition along the $\P^2$
corresponding to the intersection of lines of the configuration equals
the sum of the multiplicities of those lines.

(iii) The intersection of the point@-conditions $W'$ in $V'$ consists of
the proper transforms of the distinguished $\P^5$'s.\endproclaim

\demo{Proof} (i) is clear, while (ii) and (iii) require a computation.
Statement (iii) is vacuously true in the complement of the exceptional
divisors of the first set of blow@-ups, so we only need to check it
along a component of the exceptional divisor. Consider then
a point of intersection of two lines of the configuration; without
loss of generality we may assume this is the point $(1:0:0)$. Write
the equation of $C$ as
$$\prod_i(\beta_i y+\gamma_i z)^{r_i}\cdot \prod_j(\alpha_j x+\beta_j
y+\gamma_j z)^{r_j}=0$$
with $\alpha_j\ne 0$. We can study the blow@-up over the affine piece
$(1 : p_1 : p_2 : p_3 : p_4 : p_5 : p_6 : p_7 : p_8)$;
in these coordinates, the point@-condition corresponding to $(x_0:y_0:
z_0)$ has equation
$$\multline
\prod_i(\beta_i (p_3 x_0+p_4 y_0+p_5 z_0)+\gamma_i (p_6 x_0+p_7
y_0+p_8 z_0))^{r_i}\\
\cdot \prod_j(\alpha_j(x_0+ p_1 y_0+ p_2 z_0)
+\beta_j (p_3 x_0+p_4 y_0+p_5 z_0)+\gamma_j (p_6 x_0+p_7y_0+p_8 z_0))^{r_j}
\endmultline$$
and $(1:0:0)$ corresponds to the distinguished 2@-dimensional locus
$p_3=p_4=p_5=p_6=p_7=p_8=0$.
In a representative chart of $V'$, we can choose coordinates
$q_1,\dots, q_8$ so that the blow@-up map is given by
$$(q_1,\dots,q_8) \mapsto (q_1,q_2,q_3,q_3 q_4,\dots, q_3 q_8)\quad.$$
Therefore the point@-condition pulls back to 
$$\multline
\prod_iq_3^{r_i}(\beta_i (x_0+q_4 y_0+q_5 z_0)+\gamma_i (q_6 x_0+q_7
y_0+q_8 z_0))^{r_i}\\
\cdot \prod_j(\alpha_j(x_0+ q_1 y_0+ q_2 z_0)
+\beta_j q_3 (x_0+q_4 y_0+q_5 z_0)+\gamma_j q_3 (q_6 x_0+q_7y_0+q_8 
z_0))^{r_j}
\endmultline$$
so it contains the exceptional divisor $q_3=0$ with multiplicity $\sum
r_i$, proving (ii). Also, its proper transform has intersection with
the exceptional divisor
$$\prod_i (\beta_i (x_0+q_4 y_0+q_5 z_0)+\gamma_i (q_6 x_0+q_7 y_0+q_8
z_0))^{r_i}\cdot (x_0+ q_1 y_0+ q_2 z_0)^{\sum r_j}\quad.\tag*$$
A point $(q_1,q_2,0,q_4,\dots,q_8)$ is in the intersection of all
point@-conditions in $V'$ if (*) is identically zero for all
$(x_0,y_0,z_0)$. From (*) we see that the set of such points is given
by the equations 
$$\left\{\aligned\beta_i +\gamma_i q_6 &= 0 \\
\beta_i q_4 +\gamma_i q_7 &= 0 \\
\beta_i q_5 +\gamma_i q_8 &= 0
\endaligned\right.$$
It is easy to check that these are precisely the equations of the
proper transforms of the distinguished $\P^5$'s corresponding to those
lines in the configuration which contain $(1:0:0)$, and this verifies
(iii).\qed\enddemo

\proclaim{Claim 2.3} (i) The point@-conditions in $V'$ contain the
proper transform of a distinguished $\P^5$ with multiplicity equal to
the multiplicity of the corresponding line in the configuration.

(ii) The intersection of all point@-conditions in $\Til V$
is empty.\endproclaim
\demo{Proof} $\Til V$ is obtained from $V'$ by blowing up the
proper transforms of the distinguished $\P^5$'s. The statement can be
checked by another straightforward coordinate computation, which we
leave to the reader.\qed\enddemo

This claim
implies Proposition~2.1, as $\tilde c$ is the map corresponding to the
linear system generated by the point@-conditions, and this linear
system has empty base locus by the claim.

\remark{Remark 2.4} The orbit@-closure of a star is dominated via $\Til
c$ by the proper transform $\Til\P^5$ of a general $\P^5\subset \P^8$;
such a $\P^5$ is for example the closure of the set of matrices of
rank $2$ whose image is a fixed line not containing the center of the
star. Indeed, this locus avoids the distinguished $\P^2$
corresponding to the center of the star, and intersects transversally
the other centers of the blow@-ups considered in this section.
This implies
that the class of $\Til\P^5$ in $\Til V$ is $\Til H^3$, and hence
$$\frac 1{j!}\int_{\Til\P^5} \Til H^{5-j} \Til W^j=\frac
1{j!}\int_{\Til V} \Til H^{8-j} \Til W^j\qquad (j=0,\dots,5)$$
is the coefficient of $t^j$ in the adjusted predegree polynomial for
the star, for such a $\P^5$. This fact is used in \cite{5}.
\endremark

\subheading{\S2.3 Predegree computations} For an arbitrary
configuration of lines $C$, the construction of \S2.2 provides us with
a nonsingular variety $\Til V$ and a regular map
$$\tilde c: \Til V @>>> \P^N$$
whose image is the orbit closure of $C$. Intersection theory allows us
to compute the degree of the image of this map, and more generally to
compute the `adjusted predegree polynomial' of $C$:
$$P(t)=\sum_{j\ge 0} \left(\int \Til H^{8-j} \Til W^j\right)
\frac{t^j}{j!}\quad.$$
For stars and fans, the result of this computation can be stated quite
succinctly:

\proclaim{Theorem 2.5} (i) The adjusted predegree polynomial of a star of
lines $L_i$, appearing with multiplicity $r_i$, is 
$$\left\{\prod_i \left(1+r_i t+\frac{r_i^2
t^2}2\right)\right\}_5\quad,$$
where $\{\}_5$ denotes the truncation to $t^5$.
In particular, a star with three or more lines has predegree
$$30(e_2 e_3-e_1 e_4-e_5)\quad,$$
where $e_j$ denotes the $j$@-th elementary symmetric function in the
multiplicities $r_i$ of the lines in the star.

(ii) Let $C'$ be a star of lines with multiplicities
$r_i$, and let $C$ be a fan obtained as the union of $C'$ with a
transversal line with multiplicity $r$. Then the adjusted predegree
polynomial of $C$ is 
$$\left(1+r t+\frac{r^2 t^2}2\right)\left\{\prod_i \left(1+r_i
t+\frac{r_i^2 t^2}2\right)\right\}_5\quad.$$
In particular, if the star has three or more lines, then the predegree
of the fan is
$$630 r^2 (e_2 e_3-e_1 e_4-e_5)\quad,$$
where $e_j$ denotes the $j$@-th elementary symmetric function in the
multiplicities $r_i$ of the lines in the star.\endproclaim

For example, according to (i) the adjusted predegree
polynomial for a star consisting of three reduced lines is
$$\left\{\left(1+t+\frac{t^2}2\right)^3\right\}_5=1+3t+\frac{9t^2}2+4
t^3+\frac{9 t^4}4+\frac{3 t^5}4\quad,$$
yielding a predegree of
$\frac 34\cdot 5!=90$.
According to our enumerative interpretation
this must be the number of ordered triples of concurrent lines
containing 5~general points; and indeed this is
$3\binom 52\binom 32=90.$
By the theorem,
the predegree of a star of $d\ge 3$ simple lines is
$$(d-2)\,(d-1)\,d\,(d^2+3 d-3)\quad;$$
for $d\ge 4$ this cannot be checked by simple combinatorics, as a star
of four or more lines has moduli.
The predegree of the 4@-dimensional orbit closure of a pair of lines
is 6 (twice its degree: this is the divisor of singular conics); and
the predegree for a single line is of course~1.

As another example that can be checked combinatorially, consider the
orbit closure of a {\it triangle\/} (item (4) in the classification in
\S1).
A triangle is a fan consisting of a star of two lines union a general
line, so according to (ii)
its adjusted predegree polynomial must be
$$1+3t+\frac{9t^2}2+4 t^3+\frac{9 t^4}4+\frac{3
t^5}4+\frac{t^6}8\quad.$$
In particular, the predegree of the 6@-dimensional orbit closure of a
triangle must be
$\frac 18 \cdot 6!=90$.
The orbit closure of a triangle is however the closure of the set of
{\it all\/} triangles, so this number should agree with the number of
ordered triples of general lines in the plane, containing 6 general
points, that is
$\binom 62 \binom 42 \binom 22 = 90$,
as expected.

\demo{Proof of Theorem~2.5} The statements given above are particular
cases of the computation for arbitrary configurations of lines
(Theorem~2.8). A particularly economical
way to package all the information relevant to stars and fans is the
following:
\proclaim{Lemma 2.6} Let $C'$ be a 
configuration of lines
with adjusted predegree polynomial $P(t)$, and let $C$ be the
configuration obtained by adding to $C'$ a line with multiplicity $r$
that intersects
$C'$ transversally. Then the adjusted predegree
polynomial of $C$ is the truncation to $t^8$ of
$$\left(1+r t+\frac{r^2 t^2}2\right) P(t)\quad.$$
\endproclaim
Assuming this statement, Theorem~2.5
can be proved as follows. First, observe that the adjusted predegree
polynomial for a single $r$@-multiple line is
$$1+r t+\frac{r^2 t^2}2$$
(in this case, the orbit is simply an $r$@-tuple embedding of $\P^2$
in $\P^N$); by the lemma, it follows that the adjusted predegree
polynomial for a configuration of general distinct lines with
multiplicities $r_1, r_2, \dots$ is the truncation to $t^8$ of
$$\prod_i\left(1+r_i t+\frac{r_i^2 t^2}2\right).$$
Now, the key remarks are that\roster
\item a star has stabilizer of dimension $\ge 3$, so its adjusted
predegree polynomial has degree $\le 5$; and
\item the combinatorics of the intersections of the lines in a
configuration only affects the coefficients of degree $\ge 6$ in the
adjusted predegree polynomial.\endroster

The first observation should be clear; the second follows from the
construction in \S2.2. More explicitly, the coefficients of degree
$\le 5$ in the adjusted predegree polynomial are intersection degrees 
$$\frac 1{j!} \int \Til H^{8-j} \Til W^j$$
in $\Til V$, with $j\le 5$; that is, they involve the intersection of
3 or more proper transforms $\Til H$ of hyperplanes from $\P^8$.
However, points of intersection of the lines in the configuration
correspond to 2@-dimensional centers of blow@-up in the construction
of $\Til V$; the intersection of three or more general hyperplanes
avoids these centers, so the numbers $\int \Til H^{8-j} \Til W^j$ are
unaffected by how the lines meet, for $j\le 5$.

These observations imply part (i) of Theorem~2.5:
by (1), the coefficient of $t^j$ of the adjusted predegree polynomial
of a star must be 0 for $j\ge 6$; by (2), the 
adjusted predegree polynomial for
a star must agree with the adjusted predegree polynomial of a general
configuration up to the term of degree~5.

The second part of the theorem
follows then immediately from the first, by applying Lemma~2.6 again.
This proves Theorem~2.5 (modulo Lemma~2.6).
\qed\enddemo

Lemma 2.6 seems of independent interest. It follows immediately from a 
more general multiplicativity formula (Cor.~2.11)
for transversal configurations,
proved below. This in turn will follow from
the computation of the adjusted predegree polynomial for {\it
arbitrary\/} configurations of lines (Theorem~2.8), which will be
obtained by applying intersection theory to the construction given in
\S2.2. As in \cite{2} and \cite{3}, this can be done by applying
Theorem II from \cite{1}; in the form that we need here, this amounts to
the following formula.
\proclaim{Lemma 2.7} Let $B\overset i\to\hookrightarrow V$ be nonsingular
varieties; $X$, $Y$ hypersurfaces in $V$; $\Til X$, $\Til Y$ their proper
transforms in the blow@-up of $V$ along $B$. Then
$$\multline
\int_{\Til V} {\Til X}^{\dim V-j}\cdot {\Til Y}^j=\int_V X^{\dim
V-j}\cdot Y^j\\
- \int_B \frac{(m_{B,X}[B]+i^*[X])^{\dim V-j} (m_{B,Y}[B]+i^*[Y])^j}{c(N_B V)}
\endmultline$$
where $m_{B,X}$, $m_{B,Y}$ denote the multiplicities of $X$, $Y$ along $B$.
\qed\endproclaim

This formula and the `B\'ezout' product in $\P^8$ will yield the
following result.

Assume $C$ is an arbitrary degree@-$d$ configuration of lines $L_i$,
appearing with multiplicities $r_i$ (so $d=\sum r_i$). Line $L_i$
intersects the rest of the configuration in a tuple of points $p_j$;
for each of these points $p_j\in L_i$ we let
$$\rho_{\alpha,ij}=(\sum_{L_k\ni p_j, k\ne i} r_k)^\alpha -
(\sum_{L_k\ni p_j, k\ne i} r_k^\alpha)\quad.$$
Note that all $\rho_{\alpha,ij}=0$ if $L_i$ intersects the rest of the
configuration transversally; that is, these functions measure how far
from transversal the intersection is. Next, we let
$$\rho_{\alpha,i}=\sum_{p_j\in L_i} \rho_{\alpha,ij}\quad,$$
and define the following three functions:
$$S_{6,i}=-\rho_{5,i} r_i - 5 \rho_{4,i} r_i^2 + 10 \rho_{3,i}
r_i^3 + 5 \rho_{2,i} r_i^4$$
$$S_{7,i}=6 \rho_{6,i} r_i + 29 \rho_{5,i} r_i^2 - 50
\rho_{4,i} r_i^3 - 20 \rho_{3,i} r_i^4 - \rho_{2,i} r_i^5$$
$$S_{8,i}=-21 \rho_{7,i} r_i - 99 \rho_{6,i} r_i^2 + 155
\rho_{5,i} r_i^3 + 55 \rho_{4,i} r_i^4 + \rho_{3,i} r_i^5 -
\rho_{2,i} r_i^6$$
Again, note that these expressions vanish if $L_i$ meets the rest of
the configuration transversally.
 
\proclaim{Theorem 2.8} With the above notations, the adjusted predegree
polynomial of a configurations of lines $L_i$ is{\eightpoint
$$
\prod_i\left(1+r_i t+\frac{r_i^2 t^2}2\right)
+ \sum_i \left(1+(d-r_i) t+\frac{(d-r_i)^2 t^2}2\right)
\left(\frac{S_{6,i} t^6}{6!}+\frac{S_{7,i} t^7}{7!}+
\frac{S_{8,i} t^8}{8!} \right)
$$}
(truncated to the $t^8$ term).\endproclaim

\demo{Proof} This can be viewed as an excess intersection problem in
$\P^8$. Lemma 2.7
is used to evaluate the contribution to the B\'ezout numbers due to
the excess components, supported on the base locus of $c$. The
discussion in \S2.2 identifies a contribution from the first set of
blow@-ups (due to the intersection points of the configuration), and a
contribution from the second set of blow@-ups (due to the lines of
the configuration). Here is the raw form that the answer takes by
applying Lemma 2.7:
\proclaim{Lemma 2.9} Denote by $L_i$ the lines of the configuration $C$,
and by $p_j$ the intersection points. Let $r_i$ be the multiplicity of
$L_i$ in $C$, and let $m_j$ be the multiplicity of $p_j$ in $C$. Let
$d=\sum r_i$ be the degree of $C$. Then:

(1) the contribution due to the first set of blow@-ups is
$$\sum_{p_j}\left(\frac{m_j^6 t^6}{6!} + \frac{(7dm_j^6 - 6m_j^7)t^7}{7!}
+ \frac{(28d^2m_j^6 - 48dm_j^7 + 21m_j^8)t^8}{8!}\right)$$

(2) the contribution due to the second set of blow@-ups is{\eightpoint
$$\multline
\sum_{L_i}r_i^3\left(\frac{t^3}{3!} + (4d -3r_i)\frac{t^4}{4!} +
(10d^2- 15dr_i+6r_i^2)\frac{t^5}{5!}\right.+(20d^3 - 45d^2r_i \\
+ 36dr_i^2-10 r_i^3)\frac{t^6}{6!}
+ (35d^4 - 105d^3r_i + 126d^2r_i^2-70 d r_i^3+15
r_i^4)\frac{t^7}{7!}\\
\left. + (56d^5 - 210d^4r_i + 336d^3r_i^2-280 d^2
r_i^3+120 d r_i^4-21 r_i^5)\frac{t^8}{8!}\right)
\endmultline$$}
plus{\eightpoint
$$\multline
-\sum_{L_i}\sum_{p_j\in L_i}r_i^3\left(
(20m_j^3 - 45m_j^2r_i + 36m_jr_i^2 - 10 r_i^3)\frac{t^6}{6!}
+(140dm_j^3 - 105m_j^4 - 315dm_j^2r_i\right.\\
+ 210m_j^3r_i + 252dm_jr_i^2 - 126m_j^2r_i^2 - 70dr_i^3 +
15r_i^4)\frac{t^7}{7!} +(560d^2m_j^3 - 840dm_j^4\\
+ 336m_j^5 - 1260d^2m_j^2r_i + 1680dm_j^3r_i - 630m_j^4r_i +
1008d^2m_jr_i^2 - 1008dm_j^2r_i^2\\
\left. + 336m_j^3r_i^2 - 280d^2r_i^3 + 120dr_i^4 - 21
r_i^5)\frac{t^8}{8!}\right)
\endmultline$$}
\endproclaim
\demo{Proof of Lemma 2.9}
Consider one point $p$ of intersection of lines in the configuration;
as we have seen in \S2.2, $p$ contributes a $\P^2$ to the first set of
blow@-ups. Denote by $k$ the hyperplane class in this $\P^2$. Then the
contribution of $p$ to the first correction term given by applying the
formula in Lemma 2.7
is
$$\sum_\alpha\left(\int_{\P^2}\frac{k^{8-\alpha}
(m+d k)^{\alpha}}{(1+k)^6}\right)\frac {t^\alpha}{\alpha!}$$
where $m$ is the multiplicity of $p$ in $C$. Indeed, a general
hyperplane in $\P^8$ does not contain $B=\P^2$ (so its multiplicity
along $B$ is 0), while we have checked in Claim 2.2
that point@-conditions contain $B$ with the stated multiplicity;
$(1+k)^6$ is the total Chern class of the normal bundle of $\P^2$ in
$\P^8$. Evaluating the expression and summing over the points $p_j$
gives (1).

Next, consider a line $L$ in the configuration, appearing with
multiplicity $r$. In $V'$ we find the proper transform $\Til \P^5$ of
the distinguished $\P^5$ corresponding to $L$. Let $\ell$ denote the
hyperplane class in $\P^5$, and its pull@-back to $\Til \P^5$; also,
let $e_1,\dots, e_s$ denote the pull@-back to $\Til \P^5$ of the
exceptional divisors from the first set of blow@-ups, corresponding to
the intersection of $L$ with the rest of the configuration. Note that
the other exceptional divisors do not meet this $\Til \P^5$.
The following claim can be checked by using standard intersection theory:

\proclaim{Claim 2.10} (i) $c(N_{\Til \P^5} V')=(1+\ell-e_1-\dots-e_s)^3$.

(ii) Let $m_j$ denote the multiplicity of $C$ at the $j$@-th point of
intersection of $L$ with the rest of the configuration; then the
pull@-back of the class $W'$ of the point@-conditions in $V'$ to $\Til
\P^5$ is $d\ell -m_1 e_1-\dots-m_s e_s$.

(iii) For each $j$,
$$e_j^0 \mapsto 1\, ,\,e_j^1\mapsto 0\, ,\,e_j^2\mapsto 0\,
,\,e_j^3\mapsto \ell^3\, ,\, e_j^4\mapsto 3\ell^4\, ,\, e_j^5\mapsto 6
\ell^5$$
via the push@-forward from $\Til \P^5$ to $\P^5$, while $e_ie_j=0$ for
$i\ne j$.
\qed
\endproclaim

Using this information and Claim 2.3(ii),
Lemma 2.7
evaluates the contribution of $L$ to the correction term for the
second blow@-up. The individual intersection products are corrected by
$$\int_{\Til{\P}^5}\frac{\ell^{8-j}(r+d\ell -m_1 e_1-\dots-m_s e_s)^j}
{(1+\ell-e_1-\dots-e_s)^3}$$
This degree can be computed by pushing forward to $\P^5$ and using
Claim 2.10(iii).
Summing over the lines and the intersection points on each line gives
(2).
This finishes the proof of Lemma 2.9.
\qed\enddemo

The statement of Theorem 2.8
simply packages the information assembled in Lemma 2.9
in a more legible format.
This finishes the proof of Theorem 2.8.
\qed\enddemo

For example, if all lines meet transversally, then the adjusted
predegree polynomial of $C$ is the truncation to $t^8$ of
$$\prod_i\left(1+r_i t+\frac{r_i^2 t^2}2\right)\quad:$$
indeed, as observed above, in this case the terms $S_{k,i}$ all
vanish.

In particular, the predegree of a configuration of four or more simple
lines in general position is
$2520(e_4^2-2 e_2 e_6+2 e_8)$,
where $e_j$ denotes the $j$@-th elementary symmetric function in the
multiplicities $r_i$ of the lines in the configuration.
For example, a configuration of four general (simple) lines has
predegree~2520. This is in agreement with the naive combinatorial
count: the number of (ordered) 4@-tuples of lines through 8 general
points is
$\binom 82 \binom 62 \binom 42 \binom 22 = 2520$.
For five or more lines, this combinatorial argument cannot be applied
since the configuration has moduli.
If $d\ge 4$ lines have multiplicity one and intersect transversally,
then the predegree of the configuration is
$$\multline
2520\left(\binom d4^2-2 \binom d2 \binom d6+2 \binom d8\right)\\
=(d-3)\,(d-2)\,(d-1)\,d\,(d^4+6 d^3-31 d^2-36 d+105)\quad.
\endmultline$$
This is 0 for $d<4$, since the configuration then has small orbit.

Lemma~2.6
amounts to a `multiplicativity of adjusted predegree polynomials' when
one of the lines meets the rest of the configuration
transversally. In fact such a multiplicativity holds more generally
whenever two configurations $C'$, $C''$ are {\it transversal\/} in the
sense that every line of each configuration meets the other
configuration transversally:
\proclaim{Corollary 2.11} Let $C$ be the union of two configurations $C'$
and $C''$ meeting transversally, and let $P_{C'}(t)$, $P_{C''}(t)$
respectively be the adjusted predegree polynomials of the two
configurations. Then the adjusted predegree polynomial of $C$ is
$$\left\{P_{C'}(t)\cdot P_{C''}(t)\right\}_8\quad.$$
\endproclaim
\demo{Proof} Denote by $d,r_i,S$, $d',r'_j,S'$, and $d'',r''_k,S''$
the data for $C$, $C'$, $C''$ used in Theorem~2.8.
First, we observe that{\eightpoint
$$\left\{\prod_j\left(1+r'_j t+\frac{{r'_j}^2 t^2}2\right)\right\}_2
=\left(1+d' t+\frac{{d'}^2 t^2}2\right)\quad,$$}
and similarly for $C''$. Taking the product of the expressions for
$P_{C'}(t)$, $P_{C''}(t)$ given by the theorem and truncating to $t^8$
agrees then with the truncation of{\eightpoint
$$\multline
\prod_j\left(1+r'_j t+\frac{{r'_j}^2 t^2}2\right)\prod_k\left(1+r''_k
t+\frac{{r''_k}^2 t^2}2\right)\\
+ \sum_j \left(1+(d''+d'-r'_j) t+\frac{(d''+d'-r'_j)^2 t^2}2\right)
\left(\frac{S'_{6,j} t^6}{6!}+\frac{S'_{7,j} t^7}{7!}+
\frac{S'_{8,j} t^8}{8!} \right)\\
+ \sum_k \left(1+(d'+d''-r''_k) t+\frac{(d'+d''-r''_k)^2 t^2}2\right)
\left(\frac{S''_{6,k} t^6}{6!}+\frac{S''_{7,k} t^7}{7!}+
\frac{S''_{8,k} t^8}{8!} \right)
\endmultline$$}
that is, with the truncation of{\eightpoint
$$\multline
\prod_i\left(1+r_i t+\frac{r_i^2 t^2}2\right)
+ \sum_j \left(1+(d-r'_j) t+\frac{(d-r'_j)^2 t^2}2\right)
\left(\frac{S'_{6,j} t^6}{6!}+\frac{S'_{7,j} t^7}{7!}+
\frac{S'_{8,j} t^8}{8!} \right)\\
+ \sum_k \left(1+(d-r''_k) t+\frac{(d-r''_k)^2 t^2}2\right)
\left(\frac{S''_{6,k} t^6}{6!}+\frac{S''_{7,k} t^7}{7!}+
\frac{S''_{8,k} t^8}{8!} \right)\quad.
\endmultline$$}
{\tenpoint} 
Now the crucial observation is that transversal intersections do not
affect the $S$@-functions. Since every line $\ell'$ of $C'$ intersects
the whole configuration $C''$ transversally, the $S$@-function for
$\ell'$ (viewed as a line of $C$) agrees with its
$S'$@-function. Similarly, $S$@-functions for lines of $C''$ must agree
with the corresponding $S''$@-functions. Therefore, the expression
given above agrees with the truncation of{\eightpoint
$$
\prod_i\left(1+r_i t+\frac{r_i^2 t^2}2\right)
+ \sum_i \left(1+(d-r_i) t+\frac{(d-r_i)^2 t^2}2\right)
\left(\frac{S_{6,i} t^6}{6!}+\frac{S_{7,i} t^7}{7!}+
\frac{S_{8,i} t^8}{8!} \right)\quad,
$$}
that is, with the adjusted predegree polynomial for $C$.\qed\enddemo

Lemma~2.6 follows immediately from this corollary
by taking $C''$ to consist of a single line. The corollary
goes however a little further: for example, it shows that
multiplicativity holds for unions of stars, so long as each line of
the configuration belongs to only one star. 

\subheading{\S2.4 Concluding remarks}  
The `multiplicativity lemma' and other
qualitative considerations hold in fact in all dimensions.
The result is that the orbit closure of an $r$@-fold hyperplane in
$\P^n$ has adjusted predegree polynomial
$$P_n(t)=\sum_{k=0}^n \frac{r^k t^k}{k!}\quad,$$
(that is, `$\lim_{n\to\infty} P_n(t)=\exp(r t)$') and that the
adjusted predegree polynomial of a union of transversal hyperplanes is
the (truncated) product of such terms. For $n=1$, this recovers
\cite{4}, Proposition~1.3;
for $n>2$, this observation is apparently new. 
For example, it follows that the (pre)degree of the orbit closure of
an arrangement of $d\ge 5$ (reduced) planes in general position in
$\P^3$ is{\eightpoint
$$\multline
(d-4)(d-3)(d-2)(d-1)d\left({d^{10}}  + 10\,{d^9} + 65\,{d^8} - 
1015\,{d^7} + 63\,{d^6}\right.\\
\left. - 10885\,{d^5} + 190560\,{d^4} - 658885\,{d^3} + 1358936\,{d^2}
- 3034850\,d  + 3503500\right)\quad;
\endmultline$$}
\noindent as should be expected, this number is 0 for $d\le 4$ since
the orbit is small in that case. Of course, for $d=5$ this agrees with
the naive combinatorial count
$\binom{15}3\binom {12}3 \binom 93 \binom 63 \binom 33 =
168168000$.

Another issue of some interest is the study of linear orbits of curves
in positive characteristic. 
The classification of small orbits is more difficult in this case,
since there are many more additive polynomials. As an example, the
(non-diagonalized) multiplicative automorphism group of the
cuspidal cubic $x^2z+y^3+y^2x$ in characteristic $0$ reduces to the
additive group in characteristic $3$.
Another phenomenon appearing in positive characteristic is the 
presence of curves whose general point is an inflection point; this 
may affect the dimension of the stabilizer of a curve. For example, 
let $p>0$ be the characteristic of the ground field, and $n=p^r$ 
for $r\ge 1$; then on the curve
$x^n=y z^{n-1}$
every nonsingular point is an $n$@-flex. (In fact the curve is
{\it strange,\/} since all tangent lines contain the point $(1:0:0)$.)
The point
$(0:0:1)$, which in other characteristics is distinguished by being an
$n$@-flex, becomes indistinguishable from any other nonsingular point.
In practice this enlarges the stabilizer of the curve, making its
orbit `even smaller'. The intersection theoretic computations of
\cite{3} can be carried out for these curves, yielding their
adjusted predegree polynomial:
$$1+n t+\frac{n^2 t^2}2+\frac{n^3 t^3}{3!}+\frac{n^4 t^4}{4!}
+\frac{n^3(n^2-3) t^5}{5!}+\frac{3 n^3 (n-1) (n-2) t^6}{6!}\quad.$$
An additional subtlety is that the stabilizer for this curve is {\it
nonreduced;\/} taking this into account, the degree of the orbit
closure of the curve is then
$3(n-1)(n-2)=6\binom{n-1}2.$
For example,
the degree of the orbit closure of the 
`strange cuspidal curve' $x^3=y z^2$ in characteristic~3
equals 6,
in agreement
with a computation in \cite{6}.

\enddocument